\documentclass[10pt]{amsart}
\usepackage{amssymb,latexsym}
\theoremstyle{plain}
\newtheorem{theorem}{Theorem}

\newtheorem*{TA}{Theorem A}
\newtheorem*{TB}{Theorem B}

\newtheorem*{TD}{Theorem D}
\newtheorem*{TE}{Theorem E}
\newtheorem*{TF}{Theorem F}
\newtheorem*{TG}{Theorem G}
\newtheorem*{TI}{Theorem I}
\newtheorem{lemma}{Lemma}

\theoremstyle{definition}
\newtheorem{definition}{Definition}

\theoremstyle{remark}

\numberwithin{equation}{section}
\newcommand{\R}{\mathbb R}

\begin{document}
\title[On the persistence properties of solutions of dispersive equations]{
On the persistence properties of solutions of nonlinear dispersive equations  in weighted Sobolev spaces }
%%%%%%%%%%%%%%%%%%
%Author information
\author{J. Nahas}
\address[J. Nahas]{\'Ecole Polytechnique F\'ed\'erale de Lausanne\\
MA B1 487\\ 
CH-1015 Lausanne, Switzerland}
\email{joules.nahas@epfl.ch}
%\thanks{The second author is supported  by an  NSF grant.}
%%%%%%%%%%%%%%%%%%%%
\author{G. Ponce}
\address[G. Ponce]{Department of Mathematics\\ South Hall, Room 6607\\ University of California,
Santa Barbara\\CA 93106, USA.}
\email{ponce@math.ucsb.edu}
%\thanks{}
%%%%%%%%%%%%%
%\keywords{Parabolic}
%\subjclass{Primary: subject; Secondary: subject}
%\date{}
%\dedicatory{}
%%%%%%%%%%%%%%
\begin{abstract} We study persistence properties of solutions to some canonical dispersive models, namely the semi-linear Schr\"odinger
equation, the $k$-generalized Korteweg-de Vries equation and the Benjamin-Ono equation, in  weighted Sobolev spaces
$H^s(\R^n)\cap L^2(|x|^ldx),\;s,\,l>0$. \end{abstract}
\maketitle

\section{Introduction}\label{S:1}

This work is concerned  with  persistence properties of solutions to some nonlinear dispersive equations in weighted Sobolev spaces
$\,H^s(\R^n)\cap L^2(|x|^ldx), s, l>0$.
We shall consider the initial value problems (IVP) associated to  the following dispersive models : the nonlinear Schr\"odinger (NLS) equation
\begin{equation}
\label{NLS}
i \partial_tu +\Delta u=\mu|u|^{a-1}u,\;\;\;\;t\in\R,\;\;\;x\in\R^n,\;\;\;\mu=\pm 1,\,\;\;\;a>1,
\end{equation}
the $k$-generalized Korteweg-de Vries ($k$-gKdV) equations 
\begin{equation}
\label{kgKdV}
\partial_tu +\partial_x^3 u+ u^k\partial_xu=0,\;\;\;\;t,\,x\in\R,\;\;k\in\mathbb Z^+,
\end{equation}
and the Benjamin-Ono (BO) equation
\begin{equation}
\label{BO}
\partial_t u + \mathcal  H\partial_x^2u +u\partial_x u = 0, \qquad t, \,x
\in \R,
\end{equation}
where $ \mathcal  H$ denotes  the Hilbert transform
\begin{equation}
\label{hilbertt}
\mathcal H f(x)=\frac{1}{\pi}\lim_{\epsilon\downarrow 0}\int_{|y|\geq \epsilon}\frac{f(x-y)}{y}dy=-i\,(\text{sgn}(\xi)\,\widehat {f}(\xi))^{\lor}(x).
\end{equation}

These models have been widely studied in several contexts. For example, the KdV $k=1$ in \eqref{kgKdV}
was first deduced as a model for long waves propagating in a channel. Subsequently the KdV and its modified form
($k=2$ in \eqref{kgKdV}) were found to be relevant in a number of different
physical systems. Also they have been studied because of
their relation to inverse scattering theory \cite{GGKM}. The NLS arises as a model in several different physical phenomena (see \cite{SuSu} and references therein). In the particular,  case 
$n=1$ and $a=3$ it has been shown to be completely integrable \cite{ZS}. The BO equation \eqref{BO} was first deduced  in \cite{Be} and  \cite{On} as  a model for long internal gravity waves in deep stratified fluids.
It was also shown that it is  a completely integrable system (see \cite{AbFo}, \cite{CoWi} and references therein).

 We recall the notion of  well posedness given  in \cite{Ka1} : the  IVP  is said to be locally well posed (LWP) in
the function space $X$ if for each $u_0\in X$ there exist $T>0$ and a unique solution
$u\in C([-T,T]:X)\cap....=Y_T$ of the equation, with the map data $\to$ solution being locally continuous
from $X$ to $ Y_T$. This notion of LWP includes the \lq\lq  persistent" property, i.e. the solution
describes a continuous curve on $X$. In particular, it implies that the solution flow defines a dynamical
system in $X$. When $T$ can be taken arbitrarily large  one says that the corresponding IVP is globally well posed (GWP) in $X$.

First, we shall study the Schr\"odinger equation \eqref{NLS}.
\section{The Schr\"odinger equation \eqref{NLS}}\label{S:2}
\vskip.1in

The results in \cite{CW1}, \cite{CW2}, \cite{GV1}, \cite{Ka2}, and \cite{Ts} yield the following LWP
theory in the classical Sobolev spaces
$H^s(\R^n)$ for the IVP associated to the NLS equation \eqref{NLS}.
\vskip.1in 
\begin{TA}\label{Theorem 1} Let $s_c=n/2-2/(a-1)$.
\begin{enumerate}
\item[(I)] If $s>s_c $, $s\geq 0$, with $[s]\leq a-1$ if $a$ is not an odd integer, then for each $u_0\in H^s(\R^n)$
there exist $T=T(\|u_0\|_{s,2})>0$ and a unique solution $u=u(x,t) $ of the IVP associated to the NLS equation \eqref{NLS}
with
\begin{equation}
\label{1}
u\in C([-T,T]:H^s(\R^n))\cap L^q([-T,T]:L^p_s(\R^n)) =Z^{s}_T.
\end{equation}
Moreover, the map data $\to$ solution is locally continuous from $H^s(\R^n)$ into  $Z^{s}_T$.

\item[(II)] If $s=s_c$ and $s\geq 0$, then part (I) holds with $T=T(u_0)>0$.
\end{enumerate}

\end{TA}
\vskip.1in

  \underbar{Notations} : (a) for  $1<p<\infty$ and $s\in \R$
  \begin{equation}
 \label{ps}
 L^p_s(\R^n)\equiv (1-\Delta)^{-s/2} L^p(\R^n)=J^{-s/2} L^p(\R^n),\;\;\;\;\;\|\cdot\|_{s,p}\equiv \|(1-\Delta)^{s}\cdot\|_p,
 \end{equation}
 with $L^2_s(\R^n)=H^s(\R^n)$,

 (b) the pair of indices $(q,p)$ in \eqref{1} are given by the Strichartz estimates (see  \cite{Str2} and  \cite{GV1}):
 \begin{equation}
 \label{str}
( \int_{-\infty}^{\infty} \|e^{it\Delta}u_0\|_p^qdt)^{1/q}\leq c\|u_0\|_2,
 \end{equation}
 where
 $$
 \frac{n}{2}=\frac{2}{q}+\frac{n}{p},\;\;\;\;\;2\leq p\leq \infty,\;\;\text{if}\;\;\, n=1,\;\;2
 \leq p<2n/(n-2),\;\;\;\text{if}\;\;\,n\geq 2.
 $$

\vskip.1in
 The value  $s_c=n/2-2/(a-1)$ in Theorem A is determined  by a scaling argument  : if $u(x,t) $ is a solution of the 
 IVP  associated to the NLS equation \eqref{NLS}, then $u_{\lambda}(x,t)=\lambda^{2/(a-1)}u(\lambda x,\lambda^2t)$ satisfies the same equation
 with data $u_{\lambda}(x,0)=\lambda^{2/(a-1)}u_0(\lambda x)$. Hence, for $s\in\R$
 \begin{equation}
 \label{scaling}
 \|D^{s}u_{\lambda}(x,0)\|_2 =c \| |\xi|^s \widehat {u_{\lambda}}(\xi,0)\|_2=c\lambda^{2/(a-1)+s-n/2}\|u_0\|_2,
 \end{equation}
 is independent of $\lambda $ when $s=s_c$. In Theorem A the case (I) corresponds to the
 sub-critical case and (II) to the critical one.  In the latter, one has that if $\|D^{s_c}u_0\|_2$ is sufficiently small, then the local solution extends globally in time. 

 For the optimality of the results in Theorem A see 
\cite{BKPSV},  \cite{CCT}, and \cite{KPV2}.

Formally, solutions of the NLS equation \eqref{NLS} satisfies the following conservation laws:
$$
\|u(\cdot, t)\|_2=\|u_0\|_2,
$$
and
$$
E(t)=\int_{\R^n}( |\nabla_xu(x,t)|^2 +\frac{2\mu}{a+1} |u(x,t)|^{a+1})dx=E(0).
$$
Using these conservation laws one can extend the LWP results in Theorem A to a GWP one, for details we refer to \cite{Bo2}, \cite{Ta3}, and references therein.

 Concerning the persistence properties in weighted Sobolev spaces  of solutions of the IVP associated to the NLS equation \eqref{NLS}  one has the following result established in  \cite{HNT1}, \cite{HNT2}, and \cite{HNT3}.

  \vskip.1in
\begin{TB}
\label{Theorem 2}
 In addition to the hypothesis in Theorem A assume   $u_0\in L^2(|x|^{2m}dx)$, $m\in \mathbb Z^+$
 with $m\leq  a-1$ if $a$ is not an odd integer.
\begin{enumerate}
\item[(I)] If $\;s\geq m$, then
 \begin{equation}
\label{4a}
u\in C([-T,T]:H^s\cap L^2(|x|^{2m}dx))\cap L^q([-T,T]:L^p_s\cap L^p(|x|^{2m}dx)=Z_T^{s,m}.
\end{equation}

\item[(II)] If $\;1\leq s< m$, then \eqref{4a} holds with $[s]$ instead of $m$ and
 \begin{equation}
\label{5a}
\Gamma^{\beta}   u=(x_j+2it\partial_{x_j})^{\beta}     u\in C([-T,T]:L^2)\cap L^q([-T,T]:L^p),
\end{equation}
for any $\beta \in (\mathbb Z^+)^n$ with $|\beta|\leq m$.
\end{enumerate}
 \end{TB}

 \vskip.1in
The proof of Theorem B (see  \cite{HNT1}, \cite{HNT2}, \cite{HNT3}) combines the operators (\lq\lq vector fields")
\begin{equation}
\label{identity}
\Gamma_j =x_j+2it\partial_{x_j} =e^{i|x|^2/4t} 2it\partial_{x_j}(e^{-i|x|^2/4t}\cdot)=
e^{it\Delta}x_j e^{-it\Delta}\cdot\,,\;\;\;j=1,..,n,
\end{equation}
their commutative relation
\begin{equation}
\label{identity2}
(i \partial_t +\Delta)\Gamma_j u=\Gamma_j(i \partial_tu +\Delta u),\;\;\;\;\;\;j=1,..,n,
\end{equation}
so that $ \,e^{it\Delta}(x_ju_0)=\Gamma_j\,e^{it\Delta}u_0$, and the
structure of the nonlinearity in \eqref{NLS}.

It should be remarked that Theorem B shows that the amount of decay in  $L^2(|x|^{2m}dx)$ preserved by the solution depends on the regularity in the Sobolev scale  $H^s,\,s\geq 0$) of
the data, and the non-preserved decay is transformed in \lq\lq local regularity". In particular, \eqref{5a}
tells us that  $\,t^{\beta} \partial_x^{\beta}u\in L^2_{loc}(\mathbb R^n)$, for $ |\beta|\leq m$ and $t\in[-T,T]-\{0\}$.

Also one notices  that the power of the weight $m$ in Theorem B  is assumed to be an
integer. In  \cite{NP} we  were able to  remove this restriction. 

\vskip.1in
\begin{theorem}
\label{Theorem 3}
 In addition to the hypothesis in Theorem A  assume  $u_0\in L^2(|x|^{2m}dx)$, $m>0$
 with $[m] \leq a-1$ if $a$ is not an odd integer.
\begin{enumerate}
\item[(I)] If $s\geq m$,
 \begin{equation}
\label{4}
u\in C([-T,T]:H^s\cap L^2(|x|^{2m}dx))\cap L^q([-T,T]:L^p_s\cap L^p(|x|^{2m}dx)=Z_T^{s,m}.
\end{equation}

\item[(II)] If $1\leq s< m$, then \eqref{4} holds with $[s]$ instead of $m$ and
 \begin{equation}
\label{6}
\varGamma^b\,\Gamma^{\beta}   u(\cdot,t)\in  C([-T,T]:L^2)\cap L^q([-T,T]:L^p),
\end{equation}
where $ \varGamma^b   =  e^{i|x|^2/4t} 2^b t^b    D^b    \left(e^{-i|x|^2/4t}\,\cdot\right)$ with $|\beta|=[m]$ and $b=m-[m]$. 

In particular,
 \begin{equation}
\label{6b}
t^{m}\,\partial_x^{\beta}\,D^{b}    u(\cdot,t)\in L^2_{loc}(\R^n),\;\;\;\;|\beta|=[m],\;\;b=m-[m],\;\;t\in(-T,T)-\{0\}.
\end{equation}
\end{enumerate}
 \end{theorem}
\vskip.1in 
As an application of this result we also prove that the persistence property in these weighted spaces can only hold for regular enough solutions. More precisely: 
  \begin{lemma}
 \label{lemma5}
Let $u$ be a solution of the IVP associated to the NLS equation \eqref{NLS}  provided by Theorem A.
If there exist two times $t_1, t_2\in[0,T]$, $t_1\neq t_2$ such that
\begin{equation}
\label{condi}
|x|^mu(t_1),\;\;\;|x|^mu(t_2)\in L^2(\mathbb R^n),\;\;\;\;\;\;\;m>s,
\end{equation}
$m\leq a-1$ if $a$ is not an odd integer, then
$$
u\in C([-T,T]:H^m\cap L^2(|x|^{2m}dx))\cap L^q([-T,T]:L^p_m\cap L^p(|x|^{2m}dx).
$$

Moreover, if $a$ is an odd integer and \eqref{condi} holds for all $m\in\mathbb Z^{+}$, then
\begin{equation}
\label{sss}
u\in  C([-T,T]: \mathbb S(\R^n)).
\end{equation}
 \end{lemma}

\vskip.1in 
 A  key ingredient in our proof  was an appropriate version of the Leibnitz rule for
 homogeneous fractional derivatives of order $b\in\R$ 
 \begin{equation}
 \label{derc}
 D^b   f(x)\equiv ((2\pi |\xi|)^b    \hat f)^{\lor}(x)
 \end{equation}
deduced as a direct consequence of  the  characterization of the $L^p_s(\R^n)$ spaces (see \eqref{ps})
 given  in \cite{St1}.

 \begin{TD}
 \label{Theorem 4}
Let $b\in (0,1)$ and $\; 2n/(n+2b)\leq p< \infty$. Then $f\in  L^p_b(\R^n)$ if and only if
\begin{equation}
\label{d1}
\begin{aligned}
&\;(a)\;\, f\in L^p(\R^n),\\
\\
&\;(b)\;\;\;\;\;\mathcal D^b f(x)=(\int_{\R^n}\frac{|f(x)-f(y)|^2}{|x-y|^{n+2b   }}dy)^{1/2}\in L^p(\R^n),
\end{aligned}
\end{equation}
with
 \begin{equation}
\label{d1-norm}
\|f\|_{b,p}= \|(1-\Delta)^{b/2} f\|_p\simeq \|f\|_p+\|D^b    f\|_p\simeq \|f\|_p+\|\mathcal D^b       f\|_p.
\end{equation}
 \end{TD}
\vskip.1in 

  For the proof of Theorem D we refer to \cite{St1}, where the optimality of the lower bound $2n/(n+2b)$ 
  was also established. The case $p=2n/(n+2b)$ was proven in \cite{F}. 
  For  a detailed discussion on the  different characterizations of the $L^p_s(\R^n)$ spaces we refer to
 \cite{St1} and  \cite{Str1}. 

 It is easy to see that  for $p=2$ and $b\in(0,1)$ one has
\begin{equation}
\label{pointwise3}
\|  \mathcal D^b f\|_2\simeq \|D^b      f\|_2,
\end{equation}
\begin{equation}
\label{pointwise2}
\|  \mathcal D^b      (fg)\|_2\leq c(\|f\, \mathcal D^b      g\|_2  +\|g\, \mathcal D^b      f\|_2),
\end{equation}
and for $p> 2n/(n+2b)$  \begin{equation}
 \label{a0}
    \mathcal D^b    (fg)(x) \leq \|f\|_{\infty} \,\mathcal D^b    g(x)  +|g(x)|\, \mathcal D^b      f(x).
  \end{equation}
  
We observe that in \eqref{pointwise2} both terms on the right hand side are  estimates on the  product
 of functions.
 We do not know whether or not \eqref{pointwise2} still holds with $D^b$ instead of $\mathcal D^b\,$,or for $\,p\ne 2$ .
 
Theorem D (i.e. the estimates \eqref{pointwise2}-\eqref{pointwise3}) allows us to get  the following inequalities:

--(i)  Let $b \in (0,1)$. For  any $t>0$
 \begin{equation}
 \label{aa1}
  \mathcal D^b      (e^{it|x|^2})\leq c(t^{b/2}+t^b    |x|^b   ).
  \end{equation}
  
--(ii)  Let $b\in(0,1)$. Then there exists $c=c(b)>0$ such that for  any $t\in \R$
 \begin{equation}
 \label{aa2}
  \| |x|^b    \,e^{it\Delta} f\|_2\leq c (t^{b/2}\|f\|_2+ t^b   \|D^b   f\|_2+ \| |x|^b    f\|_2).
  \end{equation}

--(iii) Defining   the operator $\varGamma^b $  for $b>0$ as in Theorem \ref{Theorem 3} (see \eqref{6})
\begin{equation}
\label{de1}
\varGamma^b   \equiv  \varGamma^b   (t)= e^{i|x|^2/4t} 2^b t^b    D^b    \left(e^{-i|x|^2/4t}\,\cdot\right),
\end{equation}
one has for  $b >0$ and $t\in \R$ that
 \begin{equation}
 \label{aa3}
 \varGamma^b   (t) e^{it\Delta} f = e^{it\Delta}(|x|^b   f),
   \end{equation}
   and consequently
   \begin{equation}
 \label{aa4}
 \varGamma^b   (t)  f = e^{it\Delta}(|x|^b  e^{-it\Delta} f).
 \end{equation}

In addition to the estimates \eqref{aa1}-\eqref{aa4} the following two lemmas were 
essential in the proof of Theorem \ref{Theorem 3} given in \cite{NP}. 
The first is a version of the Gagliardo-Nirenberg inequality for fractional derivatives.
 
\begin{lemma}
\label{Lemma G-N} 
Let $1<q,p, r <\infty$ and $\;0<\alpha<\beta$. Then 
\begin{equation}
\label{gn1}
\| D^{\alpha} f\|_p\leq c\|f\|_r^{1-\theta}\,\|D^{\beta} f\|_q^{\theta},
\end{equation}
with
\begin{equation}
\label{gn2}
\frac{1}{p}-\frac{\alpha}{n}=(1-\theta) \frac{1}{r} +\theta\left(\frac{1}{q}-\frac{\beta}{n}\right),\;\;\;\;\;\;\;\;\theta\in[\alpha/\beta,1].
\end{equation}
\end{lemma}

The second is an interpolation estimate, which as Lemma \ref{Lemma G-N},  is  a consequence of the three line theorem.

\begin{lemma}

\label{lemma1}
Let $a,\,b>0$. Assume that $ J^af=(1-\Delta)^{a/2}f\in L^2(\mathbb R)$ and \newline $\langle x\rangle^{b}f=
(1+|x|^2)^{b/2}f\in L^2(\mathbb R)$. Then for any $\theta \in (0,1)$
\begin{equation}
\label{complex}
\|J^{ \theta a}(\langle x\rangle^{(1-\theta) b} f)\|_2\leq c\|\langle x\rangle^b f\|_2^{1-\theta}\,\|J^af\|_2^{\theta}.
\end{equation}
\end{lemma}
\vskip.1in

For the study of persistence properties of the solution to the IVP associated to the NLS equation
 \eqref{NLS} in exponential weighted spaces we refer to \cite{EKPV1}, \cite{EKPV2}, and references therein. 

\vskip.1in

Next, we shall consider the $k$-gKdV equation \eqref{kgKdV}.
\vskip.1in

\section{The $k$-generalized Korteweg-de Vries equation \eqref{kgKdV}}\label{S:3}
\vskip.1in

The following theorem describes the LWP theory in the classical Sobolev spaces $H^s(\R)$ 
for the IVP associated to the $k$gKdV equation \eqref{kgKdV}.

\begin{TE}\label{Theorem 11} 

\begin{enumerate}

\item[(I)]  The IVP associated to the equation \eqref{kgKdV} with $k=1$  is LWP in $H^s(\R)$ for $s\geq s^*_1=-3/4$.

\item[(II)]   The IVP associated to the equation \eqref{kgKdV} with $k=2$  is LWP in $H^s(\R)$ for $s\geq s^*_2=1/4$.

\item[(III)]   The IVP associated to the equation \eqref{kgKdV} with $k=3$  is LWP in $H^s(\R)$ for $s\geq s^*_3=-1/6$.

\item[(IV)]   The IVP associated to the equation \eqref{kgKdV} with $k\geq 4$  is LWP in $H^s(\R)$ for $s\geq s^*_k=(k-4)/2k$.
\end{enumerate}
\end{TE}

The result $s>-3/4$ for the case $k=1$ was established in \cite{KPV3}. The limiting value $s=-3/4$ 
was obtained in \cite{CCT}. The result for the case $k=2$ was proven in \cite{KPV1}. The result $s>-1/6$ for the case $k=3$ was given in \cite{AG}. The limiting value $s=-1/6$ 
was obtained in \cite{Tao2}. The proof of the cases  $k\geq 4$ was given in \cite{KPV1}.

The above local results apply to both real and complex valued functions.

The scaling argument described in \eqref{scaling} affirms that LWP should hold for $s\geq s_k=(k-4)/2k$.
As Theorem E shows this is the case for $k\geq 3$ (where for $s_k=s^*_k$ one has $T=T(u_0)$). However, in the cases $k=1$ and $k=2$ the values 
suggested by the scaling do not seem to be reachable  in the Sobolev scale, see \cite{KPV2}, and \cite{CCT}. For the sharpness of these results we refer to \cite{BKPSV}, 
\cite{KPV2}, and \cite{CCT}.

Real valued solutions of the $k$-gKdV equation \eqref{kgKdV}  formally satisfy at least three conservation laws:
$$
I_1(u)=\int_{-\infty}^{\infty}\,u(x,t)dx,\;\;\;\;\;\;\;\;I_2(u)=\int_{-\infty}^{\infty}\,(u(x,t))^2dx,
$$
$$
I_3(u)=\int_{-\infty}^{\infty}\,((\partial_x u(x,t))^2-\frac{2}{(k+1)(k+2)}u(x,t)^{k+2})dx.
$$
 It was proven in \cite{CKSTT} that for $k=1$ and $k=2$  one has global well posedness for $s>-3/4$ and $s>1/4$, respectively. The global cases for $k=1,\;s=-3/4$ and $k=2,\;s=1/4$ were proven in \cite{Gu} and \cite{Ki1}. 
For the case $k=3$ the global well posedness is known for $s> - 1/42$, see \cite{AGP}.

 For $k=4$ blow up of \lq\lq large" enough solutions was proven  in \cite{FMYM}. Similar results for $k\geq 5$ remain an open problem.

Concerning the persistence of these solutions in weighted Sobolev spaces one has the following result found in \cite{Ka1}.
 \vskip.1in
\begin{TF}
\label{Theorem 22}
Let $m\in\mathbb Z^+$.  Let $u\in C([-T,T]:H^{s}(\R))\cap.....$ with $s\geq 2m$ be the solution of the IVP associated to the equation \eqref{kgKdV} provided by Theorem E. If 
$u(x,0)=u_0(x)\in L^2(|x|^{2m}dx)$, then
$$
u\in C([-T,T]:H^{s}(\R)\cap L^2(|x|^{2m}dx)).
$$
 \end{TF}
 \vskip.1in
 
We recall that if  for a solution $u\in C([0,T]:H^s(\R))$ of \eqref{kgKdV} one has that $\,\exists \,t_0\in[0,T]$ such that $u(\cdot,t_0)\in H^{s'}(\R),\,s'>s$, then  $u\in C([0,T]:H^{s'}(\R))$. So we shall mainly consider the most interesting case $s=2m$ in Theorem F.

The proof of Theorem F combines the operator
$$
\Gamma =x+3t\partial_x^2,
$$
and its commutative relation with the linear part $L=\partial_t+\partial_x^3$ of the equation \eqref{kgKdV}
i.e. 
$$
\Gamma (\partial_t+\partial_x^3)v=(\partial_t+\partial_x^3)  \Gamma v.
$$

As in the case of the NLS equation \eqref{NLS} we would like to extend Theorem F  where $m\in\mathbb Z^+$
to the case $m\in\R,\, m>0$. Our first result in this direction is the following:

\vskip.1in
\begin{theorem}
\label{Theorem 23}
Let $m\geq 0$.  Let $u\in C([-T,T]:H^{m}(\R))\cap.....$ with $m\geq \max\{s^*_k;0\}$ be the solution of the IVP associated to the equation \eqref{kgKdV} provided by Theorem E. If 
$u(x,0)=u_0(x)\in L^2(|x|^{m}dx)$, then

\begin{enumerate}
\item[(I)] If $m<1$, then for any $\epsilon>0$ 
$$
u\in C([-T,T]:H^{m}(\R)\cap L^2(|x|^{m-\epsilon}dx)).
$$
\item[(II)] If $m\geq 1$, then 
$$
u\in C([-T,T]:H^{m}(\R)\cap L^2(|x|^{m}dx)).
$$
\end{enumerate}
\end{theorem}

\vskip.1in
In \cite{JN1} and \cite{JN2}  the loss of power $\epsilon>0$ in the weight  when $m<1$ was removed for the equation \eqref{kgKdV} with non-linearity $k=2,4,5,....$.
More precisely, the following  optimal result was established in \cite{JN2}:

\vskip.1in
\begin{theorem}
\label{Theorem 24}
Let $m\geq max\{s_k^*;0\}$ with $k=2,4,5,...$.  Let $u\in C([-T,T]:H^{m}(\R))\cap.....$  be the solution of the IVP associated to the equation \eqref{kgKdV} provided by Theorem E. If 
$u(x,0)=u_0(x)\in L^2(|x|^{m}dx)$, then
$$
u\in C([-T,T]:H^{m}(\R)\cap L^2(|x|^{m}dx)).
$$
\end{theorem}
\vskip.1in

It should be remarked that in the cases $k=1$ and $k=3$ the proof of the local theory in Theorem E is based on the spaces $X_{s,b}$ introduced in the context of dispersive equations  in \cite{Bo1}.
For all the other powers $k$ one has a local existence theory based on a contraction principle in a spaces defined by mixed norms of the type $L^p(\R : L^q([0,T]))$ or 
$L^q([0,T] : L^p(\R))$ (see \cite{KPV1}). This is the main difficulty in extending the optimal result in Theorem \ref{Theorem 24} to the powers $k=1$ and $k=3$ in \eqref{kgKdV}.
\vskip.2in
\underbar{Proof of Theorem \ref{Theorem 23}}
\vskip.1in

We shall sketch the  ideas in the proof of Theorem \ref{Theorem 23} and refer to \cite{JN1} and  \cite{JN2}  for the justification of the argument and further details.

Following Kato's idea  in \cite{Ka1} to establish the local smoothing effect (i.e. multiplying the equation \eqref{kgKdV} by $u(x,t) \phi(x)$, integrating the result, and using integration by parts) one formally gets the identity
\begin{equation}
\label{key1}
\frac{d\;}{dt}\int u^2\phi dx+ 3\int (\partial_xu)^2 \phi' dx -\int u^2 \phi^{(3)}dx -\frac{2}{k+2}\int u^{k+2}\phi'dx=0.
\end{equation}

Let us consider first the case $max\{s_k^*;0\}\leq m<1$.

From the local theory one has the following estimates for the solution $u=u(x,t)$ 
\begin{equation}
\label{21}
\sup_{x\in\R}\,( \int_0^T |\partial_x D^m_x u(x,t)|^2dt )^{1/2}< c_T \|J^m u_0\|_2=c_T\|u_0\|_{m,2},
\end{equation}
(the sharp form of the local smoothing effect found in \cite{KPV0}-\cite{KPV1}), and
\begin{equation}
\label{22}
\begin{aligned}
&\|D^m_xu\|_{L^2_xL^2_T}=(\int_{-\infty}^{\infty} \int_0^T |D^m_xu(x,t)|^2dtdx)^{1/2} \\
&\leq T^{1/2} \sup_{t\in[0,T]} \| D^m_x u(t)\|_2 < c_T \|D^m u_0\|_2\leq c_T\|u_0\|_{m,2}.
\end{aligned}
\end{equation}

Now, we consider the extensions of the estimates in \eqref{21}-\eqref{22} to the operators
$D_x^{1+m+iy}$ and $D_x^{m+iy},\;y\in\R$ respectively. First, in the linear case one has the estimates
\begin{equation}
\label{inequalities}
\begin{aligned}
& \| D_x^{m+1+iy}v\|_{L^{\infty}_xL^2_T}\leq c_T\| D^mv_0\|_2,\\
&\| D_x^{m+iy}v\|_{L^{2}_xL^2_T}\leq c_T\| D^mv_0\|_2,
\end{aligned}
\end{equation}
for 
\begin{equation}
\label{linear}
v(x,t)=U(t)v_0(x)=c\,\int_{-\infty}^{\infty} \,e^{ix\xi} e^{it\xi^3}\widehat v_0(\xi)d\xi.
\end{equation}
To apply the three line theorem 
we consider the function $F(z)$ defined on $\mathcal S = \{z\in\mathbb C\;:\; \Re(z)\in [0,1]\}$  
$$
F(z)=\int_{-\infty}^{\infty}\int_0^T D_x^{s(z)}v(x,t)\, \phi(x,z)\,f(t)\,dt dx,
$$
where
$$
s(z)=(1-z)(1+m)+z m,\;\;\;1/q(z)=(1-z)+z/2,\;\;\;q=2/(2-m),
$$
$$
\phi(x,z)=|g(x)|^{q/q(z)}\,\frac{g(x)}{|g(x)|},\;\;\;\;\text{with}\;\;\;\;\|g\|_{L^{2/(2-m)}_x}=\|f\|_{L^2([0,T])}=1,
$$
which is analytic on the interior of $\mathcal S$. So using that
$$
\|\phi(\cdot,0+iy)\|_1=\|\phi(\cdot,1+iy)\|_2=1,
$$
one gets that
\begin{equation}
\label{23}
\begin{aligned}
 \|\partial_xv\|_{L^{2/m}_xL^2_T}&\leq c \|D_xv\|_{L^{2/m}_xL^2_T}\\
&\leq c \,\sup_{y\in\R} \| D_x^{1+m+iy} v\|_{L^{\infty}_xL^2_T}^{1-m}\,\sup_{y\in\R}\|D^{m+iy}_xv\|_{L^2_xL^2_T}^m\leq c_T \| D^mv_0\|_{2}.
\end{aligned}
\end{equation}

Inserting the estimate \eqref{23} in the proof of the local well posedness one obtains that
\begin{equation}
\label{23a}
 \|\partial_xu\|_{L^{2/m}_xL^2_T}\leq c_T \|u_0\|_{m,2},
\end{equation}
for $u=u(x,t) $ solution of the $k$-gKdV equation \eqref{kgKdV}.

Now taking $\phi(x)=\langle x\rangle^{m-\epsilon},\;\epsilon>0$ sufficiently small  in \eqref{key1}, (we recall that $m<1$)  and integrating in the time interval $[0,T]$ one finds that 
\begin{equation}
\label{24}
\begin{aligned}
&\int_0^T\int_{-\infty}^{\infty} (\partial_xu(x,t))^2 \phi'(x) dx dt
=c \|\partial_x u \,\langle x\rangle^{\frac{m}{2}-\tfrac{1}{2}-\tfrac{\epsilon}{2}}\|^2_{L^{2}_xL^2_T}\\
& \leq c \| \langle x\rangle^{m/2-1/2-\epsilon/2}\|_{L^{2/(1-m)}_x} \|\partial_xu\|_{L^{2/m}_xL^2_T}\leq c_{m,\epsilon}\|\partial_xu\|_{L^{2/m}_xL^2_T},
\end{aligned}
\end{equation}
which combined with \eqref{23} and \eqref{key1} shows that $ \langle x\rangle^{m/2-\epsilon/2} u(\cdot,t)\in L^2(\R)$ for $t\in[0,T]$. This  basically completes the proof of the case $m<1$.
\vskip.1in

Next,  we shall consider the case $m\geq 1$.
\vskip.1in
We take in \eqref{key1} $\phi(x)= \langle x\rangle^{m}$ in \eqref{key1}, so we need to estimate the term
$$
\int_{-\infty}^{\infty} |\partial_xu(x,t)|^2 \langle x\rangle^{m-1}dx=\|\partial_x u(\cdot,t) \langle \cdot \rangle^{(m-1)/2}\|^2_{L^2_x}.
$$
Thus, combining Lemma 3 in the previous section, the preservation of the $L^2$-norm of the solution, and  Lemma \ref{lemma1}  it follows that
\begin{equation}
\label{27}
\begin{aligned}
&\|\partial_x u(\cdot,t) \langle \cdot \rangle^{(m-1)/2}\|_2\\
&
\leq \|\partial_x (u(\cdot,t) \langle \cdot \rangle^{(m-1)/2})\|_2+c\| u(\cdot,t) \langle \cdot \rangle^{(m-3)/2}\|_2\\
&
\leq \|\partial_x J^{-1} J (u(\cdot,t) \langle \cdot \rangle^{(m-1)/2})\|_2+c\| u(\cdot,t) \langle \cdot \rangle^{m/2}\|_2\\
&\leq c\| J (u(\cdot,t) \langle \cdot \rangle^{(m-1)/2})\|_2+c\| u(\cdot,t) \langle \cdot \rangle^{m/2}\|_2\\
&\leq 
c\| J^m u(\cdot,t)\|_2^{1/m}\| u(\cdot,t) \langle \cdot \rangle^{m/2} \|^{1-1/m}_2+c\| u(\cdot,t) \langle \cdot \rangle^{m/2}\|_2.
\end{aligned}
\end{equation}

Hence, inserting \eqref{27} in \eqref{key1}, using Young and Gronwall inequalities, the hypothesis $m\geq 1$,  and the fact that the $H^m$-norm of the solution
is bounded in the time interval $[0,T]$ one obtains the desired result
$$
\sup_{t\in[0,T]}\|\langle x\rangle^{m/2} u(\cdot,t)\|_{L^2}<\infty.
$$

This completes the sketch of the proof of Theorem \ref{Theorem 23}.

\vskip.1in

To finish this section concerning the $k$-gKdV equation \eqref{kgKdV} we will make some comments concerning  the proof of Theorem \ref{Theorem 24} given in \cite{JN1} and \cite{JN2}. 
One of the key element in that proof is the following commutator estimate:
 \begin{lemma}
 \label{Lemma  A1}
 Let $0<\alpha<1$ and $1<p<\infty$. Then for functions $f,g :\R\to \mathbb C$ one has that
 \begin{equation}
 \label{jn1}
 \|D^{\alpha}(fg)-fD^{\alpha}g\|_p \leq c\| Q_N(D^{\alpha}f)\|_{L^{\infty} l_N^1}\,\|g\|_2,
 \end{equation}
 \end{lemma}
where 
$$
\|Q_N(f)\|_{L^{\infty}l_N^1}\equiv \| \sum_{N\in\mathbb Z}|Q_N(f)|\,\|_{L^{\infty}},
$$
and
$$
Q_N(f)(x)=((\eta\left(\frac{\xi}{2^N}\right)+\eta\left(-\frac{\xi}{2^N}\right))\widehat f(\xi))^{\lor}(x),
$$
where $\eta \in C^{\infty}_0(\R)$ with $\,supp (\eta)\subseteq [1,2,2]$ so that
$$
\sum_{N\in\mathbb Z}\,(\eta\left(\frac{x}{2^N}\right)+\eta\left(-\frac{x}{2^N}\right))=1,\;\;\;\;\text{for}\;\;x\neq 0.
$$
\vskip.1in

In the proof of Theorem \ref{Theorem 24} for the case $k=2$ and $m=1/4$ (extremal case) given in \cite{JN1} Lemma \ref{Lemma A1} was combined with the inequality
$$
\|D_{\xi}^{1/8} Q_N \left(\frac{e^{it\xi^3}}{(1+\xi^2)^{1/8}}\right)\,\|_{L^{\infty}_{\xi}l^1_N}<\infty,
$$
to establish the main estimate in the proof.

\vskip.1in

For the study of persistence properties of the solution to the IVP associated to the $k$-gKdV equation
 \eqref{kgKdV} in exponential weighted spaces we refer to  \cite{KePoVe} and \cite{EKPV0} and references therein. 

\vskip.1in

Finally, we shall consider the BO equation  \eqref{BO}.

\section{The Benjamin-Ono equation \eqref{BO}}\label{S:4}
\vskip.1in

The LWP in the Sobolev spaces $H^s(\R)$ of the IVP associated to the BO equation \eqref{BO} has been largely considered :
in  \cite{ABFS} and 
\cite{Io1} LWP was established  for $s>3/2$, 
in \cite{Po} for $s\geq 3/2$, 
in \cite{KoTz1} for $s>5/4$, 
in \cite{KeKo} for $s>9/8$,
in \cite{Ta} for $s\geq 1$,
in \cite{BuPl} for $s>1/4$,
and 
in \cite{IoKe} LWP was proven in $H^s(\R)$ for $s\geq 0$.

Real valued solutions of the IVP \eqref{BO} satisfy infinitely many conservation laws (time invariant quantities), the first three are the following:
\begin{equation}
\begin{aligned}
\label{laws}
&\;I_1(u)=\int_{-\infty}^{\infty}u(x,t)dx,\;\;\;\;I_2(u)=\int_{-\infty}^{\infty}u^2(x,t)dx,\\
&\;I_3(u)=\int_{-\infty}^{\infty}\,(|D_x^{1/2}u|^2-\frac {u^3}{3})(x,t)dx,
\end{aligned}
\end{equation}
where $D_x=\mathcal H\,\partial_x$.

The $k$-conservation law $\,I_k$ provides an  \it a priori \rm estimate of the $L^2$-norm of the derivatives of order $(k-2)/2,\;k>2$  of the solution,  i.e. $\|D_x^{(k-2)/2}u(t)\|_2$. This  allows one to deduce GWP from LWP results.

In the  BO equation the dispersive effect is described
by a non-local operator and   is significantly weaker than that exhibited by the 
Korteweg-de Vries (KdV) equation, i.e. $k=1$ in \eqref{kgKdV}.
Indeed,  it was proven in \cite{MoSaTz}  that for any $s\in\R$ the map data-solution from $H^s(\R)$ to $C([0,T]:H^s(\R))$ is not locally $C^2$, 
and in  \cite{KoTz2} that it is not locally uniformly continuous. In particular, this  implies that no LWP results can be obtained by an argument based only on a  contraction  method. 

Consider the weighted Sobolev spaces
\begin{equation}
\label{spaceZ}
Z_{s,r}=H^s(\R)\cap L^2(|x|^{2r}dx),\;\;\text{and}\;\;\dot Z_{s,r}=\{ f\in Z_{s,r}\,:\,\widehat {f}(0)=0\}\;\;\;s,\,r\in\R.
\end{equation}
In \cite{Io1}  the following results were obtained:

\begin{TG}\label{Theorem BO1} 
\begin{enumerate} \item[(I)] The IVP associated to the BO equation \eqref{BO} is GWP in $Z_{2,2}$.

\item[(II)] If  $\,\widehat{u}_0(0)=0$, then the IVP associated to the BO equation \eqref{BO} is GWP 
in $\dot Z_{3,3}$.

\item [(III)] If  $u(x,t)$ is a solution of the IVP associated to the BO equation\eqref{BO} such that $u\in C([0,T]: Z_{4,4}) $ for arbitrary  $T>0$, then $u(x,t)\equiv 0$.

\end{enumerate}
\end{TG}
 
 We observe that the linear part of the equation in \eqref{BO} $L=\partial_t  + \mathcal  H\partial_x^2\,$
 commutes with the operator $\Gamma = x-2 t\mathcal H\partial_x$, i.e.
 $$
 [L;\Gamma]=L\Gamma-\Gamma L=0.
 $$
   Also, the solution $v(x,t)$ of the associated IVP 
 \begin{equation}
 \label{linearasso}
 v(x,t)=U(t)v_0(x)=e^{-it\mathcal H\partial_x^2}v_0(x)=\ (e^{-it\xi|\xi|}\,\widehat{v}_0)^{\lor}(x),
 \end{equation}
satisfies that  
 $v(\cdot, t)\in  L^2(|x|^{2k}dx), \,t\in[0,T]$, when $v_0\in Z_{k,k},\;k\in \mathbb Z^+$ for $k=1,2,......$ and
  $$
  \int_{-\infty}^{\infty}x^j\,v_0(x)dx=0,\;\;\;\;j=0, 1,..., k-3,\;\;\;\text{if}\;\;\;k\geq 3.
  $$

 In \cite{Io2} the unique continuation result  in  $Z_{4,4}$ in Theorem G  was improved:  
 \begin{TI}    
 \label{Theorem BO2} 
 Let $u\in C([0,T] : H^2(\R))$ be a solution of the IVP \eqref{BO}. If   there exist  three different times
 $\,t_1, t_2, t_3\in [0,T]$ such that 
 \begin{equation}
 \label{3timesw} 
 u(\cdot,t_j)\in Z_{4,4},\;\;\;\;j=1,2,3,\;\;\;\;\text{then}\;\;\;\;\;u(x,t)\equiv 0.
 \end{equation}
 
\end{TI}

\vskip.1in

As in the previous cases, the  goal was  to extend the results in Theorem G  and Theorem I from integer values to the continuum  optimal range of indices $(s,r)$.
 In this direction one finds the following results established in \cite{GFGP}:
 
 \begin{theorem}\label{Theorem BO3} 
 \item[(I)] Let $s\geq 1, \;r\in [0,s]$, and $\,r<5/2$. If $u_0\in Z_{s,r}$, then the solution $u(x,t)$ of the IVP associated to the BO equation \eqref{BO} satisfies that $ u\in C([0,\infty):Z_{s,r})$.
 
 \item[(II)] For  $s>9/8$  ($s\geq 3/2$), $\;r\in [0,s]$, and $\,r<5/2$ the IVP associated to the  BO equation\eqref{BO} is LWP (GWP resp.) in $Z_{s,r}$.

\item[(III)] If  $\,r\in [5/2,7/2)$ and $\,r\leq s$, then the IVP associated to the BO equation \eqref{BO} is GWP 
in $\dot Z_{s,r}$.

\end{theorem}

\begin{theorem}\label{Theorem BO4} 
 Let $u\in C([0,T] : Z_{2,2})$ be a solution of the IVP associated to the BO equation \eqref{BO}. If   there exist  two different times
 $\,t_1, t_2\in [0,T]$ such that 
 \begin{equation}
 \label{2timesw} 
 u(\cdot,t_j)\in Z_{5/2,5/2},\;\;j=1,2,\;\;\text{then}\;\;\;\widehat {u}_0(0)=0\,,\;\;(\text{so}\;\; u(\cdot, t)\in  \dot Z_{5/2,5/2}).
 \end{equation}
 
\end{theorem}
 
 \vskip.1in
 
 \begin{theorem}\label{Theorem BO5} 
 Let $u\in C([0,T] : \dot Z_{3,3})$ be a solution of the IVP \eqref{BO}. If   there exist  three different times
 $\,t_1, t_2, t_3\in [0,T]$ such that 
 \begin{equation}
 \label{3timeus} 
 u(\cdot,t_j)\in Z_{7/2,7/2},\;\;j=1,2,3,\;\;\text{then}\;\;\;u(x,t)\equiv 0.
 \end{equation}
 
\end{theorem}
\vskip.1in

\underline{Remarks} : Theorem \ref{Theorem BO4} and Theorem  \ref{Theorem BO5} show 
that the upper values of $r$  for the persisitence properties in  $Z_{s,r}$ and $\dot Z_{s,k}$ 
in Theorem \ref{Theorem BO3} are optimal.
We recall that if  $u\in C([0,T]:H^s(\R))$ is a solution of the BO equation  \eqref{BO} such  that $\,\exists \,t_0\in[0,T]$ for which $u(x,t_0)\in H^{s'}(\R),\,s'>s$, then  $u\in C([0,T]:H^{s'}(\R))$. So 
it suffices  to consider the most interesting case $s=r$ in \eqref{spaceZ}.

\vskip.in

 The proof of Theorems \ref{Theorem BO5} is based on weighted energy estimates and involves several inequalities concerning    the Hilbert transform $\mathcal H$.  
 
 Among them one finds  the $A_p$ condition introduced in \cite{Mu}.
 \begin{definition}\label{definition1} A non-negative function $w\in L^1_{loc}(\R)$ satisfies the $A_p$ inequality with $1<p<\infty\,$  if
 \begin{equation}
 \label{apb}
 \sup_{Q\;\text{interval}}\left(\frac{1}{|Q|}\int_Q w\right)\left(\frac{1}{|Q|}\int_Qw^{1-p'}\right)^{p-1}=c(w)<\infty,
 \end{equation}
 where $1/p+1/p'=1$.
 \end{definition}

    It was proven in \cite{MuHuWh} that this is a necessary and sufficient condition  for the Hilbert transform $\mathcal H$ to be bounded in  
  $L^p(w(x)dx)$ (see  \cite{MuHuWh}, ), i.e. $\;w\in A_p,\;1<p<\infty$ if and only  if
\begin{equation}
\label{a1}
( \int_{-\infty}^{\infty}|\mathcal Hf|^pw(x)dx)^{1/p}\leq c^*\, (\int_{-\infty}^{\infty} |f|^pw(x)dx)^{1/p},
\end{equation}
In the case $p=2$, a previous characterization of $w$ in \eqref{apb} was found in \cite{HeSz}. However, even though  the main case is for  $p=2$, the characterization \eqref{apb} will be the one used in the proof.
In particular, one has that in $\R$ 
\begin{equation}
\label{cond|x|}
|x|^{\alpha}\in A_p\;\;\Leftrightarrow\;\; \alpha\in (-1,p-1).
\end{equation}

In order to justify some of the arguments in the  proofs one need some further  continuity properties of the Hilbert transform. More precisely, the  proof  requires
 the constant $c^*$ in \eqref{a1} to depend only on  $c(w)$  the constant describing the $A_p$ condition
(see \eqref{apb}) and on  $p$. In \cite{Pe} precise bounds for the constant $c^*$ in \eqref{apb} were given
which are sharp in  the case $p=2$ and sufficient for the  purpose in \cite{GFGP}. 
 
 It will be essential in  
the  arguments in \cite{GFGP} that some commutator operators involving the Hilbert transform $\mathcal H$  are of   \lq\lq order zero".
More precisely, one shall use  the following estimate: 
 $\forall \,p\in(1,\infty),$ $l,\,m\in\mathbb  Z^+\cup\{0\},\,l+m\geq 1$  
$ \,\exists\, c=c(p;l;m)>0$ such that
\begin{equation}
\label{77}
 \| \partial_x^l[\mathcal H;\,a]\partial_x^m f\|_p\leq c \|\partial_x^{l+m} a\|_{\infty} \|f\|_p.
\end{equation}
 In the
case $l+m=1$,
\eqref{77} is Calder\'on's first commutator estimate
\cite{Ca}. The  case $l+m\geq 2$ of the estimate  \eqref{77} was proved in \cite{DaMcPo}.

\vskip.1in
 
{\underbar{ACKNOWLEDGMENT}}: J. N. was supported by the EAPSI NSF and JSPS program. G.P. was supported by NSF grant DMS-0800967. 
Part of this work was done while J. N.  was visiting Prof. Y. Tsutsumi at the  Department of Mathematics at Kyoto University whose hospitality he gratefully acknowledges.

\end{document}